\documentclass[a4paper, 
10pt
]{article}

\usepackage{subcaption}
\usepackage{pgfplots}
\usepackage{pgfplotstable}
\usepgfplotslibrary{groupplots}
\usepgfplotslibrary{fillbetween}
\usepackage{mathtools}
\usepackage{multicol}
\usepackage{comment}
\usepackage{booktabs}
\pgfplotsset{compat=1.5}
\usepackage{amssymb}
\usepackage{url}
\usepackage{bm}

\usepackage{pdfsync}

\usepackage{float}
\usepackage{tabularx}
\usepackage{enumerate}
\usepackage{array}
\usepackage{xspace}
\usepackage{tikz}
\usepackage{tikzsymbols}
\usetikzlibrary{calc,trees,positioning,arrows,chains,shapes.geometric,%
    decorations.pathreplacing,decorations.pathmorphing,shapes,%
    matrix,shapes.symbols, decorations.markings, patterns,fit,quotes,%
    arrows.meta}

\usepackage{authblk}

\usepackage{hyperref}
\usepackage[draft,inline,marginclue]{fixme}
\FXRegisterAuthor{mt}{amt}{MT}

\newtheorem{theorem}{Theorem}

\newcommand{\mupar}{\ensuremath{\boldsymbol{\mu}}}

\newcommand{\R}{\ensuremath{\mathbb{R}}}
\newcommand{\fulldim}{\ensuremath{\mathcal{N}}}
\newcommand{\reddim}{\ensuremath{N}}

\newcommand{\fullspace}{\ensuremath{\mathbb{X}^\fulldim}}
\newcommand{\redspace}{\ensuremath{\mathbb{X}^\fulldim_\reddim}}
\newcommand{\parspace}{\ensuremath{{\bf P}}}
\newcommand{\Xtrain}{\ensuremath{\mathbf{X}}}
\newcommand{\Xtest}{\ensuremath{\mathbf{\Bar{X}}}}
\newcommand{\ytrain}{\ensuremath{\mathbf{y}}}
\newcommand{\ytest}{\ensuremath{\mathbf{\Bar{y}}}}

\begin{document}

\title{Gaussian process approach within a data-driven POD framework for fluid
dynamics engineering problems}

\author[]{Giulio~Ortali\footnote{giulio.ortali@sissa.it}}
\author[]{Nicola~Demo\footnote{nicola.demo@sissa.it}}
\author[]{Gianluigi~Rozza\footnote{gianluigi.rozza@sissa.it}}

\affil{Mathematics Area, mathLab, SISSA, via Bonomea 265, I-34136
  Trieste, Italy}

\maketitle

\begin{abstract}
This work describes the implementation of a data-driven approach for the reduction of the complexity of parametrical partial differential equations (PDEs) employing Proper Orthogonal Decomposition (POD) and Gaussian Process Regression (GPR). This approach is applied initially to a literature case, the simulation of the stokes problems, and in the following to a real-world industrial problem, inside a shape optimization pipeline for a naval engineering problem.
\end{abstract}


\section{Introduction and motivations}
\label{sec:intro}

Reduced Order Modeling (ROM) offers a simplification for very complex computational
models~\cite{RozzaMalikDemoTezzeleGirfoglioStabileMola2018ECCOMAS},
keeping by high accuracy for the solution we want to compute but, at the same time,
reducing the computational cost to reach it. In a parametric context, this
means that we are able to approximate the solution for any new unforseen
parameter in a very efficient and rapid way. Especially for complex phenomena of which
numerical solutions are expensive to compute --- e.g. solve a Partial
Differential Equations (PDEs) ---, ROM makes a huge difference in terms of
the global computational load, allowing for a quick (or even real-time)
response and for a multiple queries approach, e.g. an optimization scenario,
otherwise not viable. 

In the Reduced Basis (RB)~\cite{rozza} method, a set of truth solutions for some selected
parameters are computed using a numerical solver over the full-order model.
Once these {\it snapshots} are collected, we combine them to extract the basis
which spans the low-dimensional space, thus we can exploit it by using the
basis as global test functions in a Galerkin method. Thanks to the low
dimensionality of the reduced order model, the related system results smaller
and faster to be solved.  Despite this advantage, RB meets some difficulties
for applications in industrial field: the necessity to extract the discrete
operators precludes the use of commercial software, as well as the effort
required to derive the equations onto the reduced space discourages its
exploitation.

Data-driven reduced order model are then gaining popularity to
tackle industrial problems, mainly for its easiness of application and for
the increasing quantity of available scientific data.  Within the data-driven
approach, it is possible avoiding to project equations and operators onto the
reduced space --- that can be, depending on the equations, very demanding and
complex --- and approximating the solution maninfold in the reduced space
with an interpolation/regression.
In this contribution, we coupled the Proper Orthogonal Decomposition (POD), a
well-known technique for computation of reduced space basis~\cite{ballarin2015supremizer}, with the
Gaussian Process Regression (GPR) for a probabilistic prediction of the new
solutions. While the POD technique is pretty standard for projection based
snapshots method in several fields of application --- e.g. structural
engineering~\cite{pichi2019reduced},
fluid dynamics~\cite{StaRo2018, hijazi2020data, StabileZancanaroRozza2020} ---, the GPR has been adopted
only in~\cite{guo2018reduced} for a structural problem. The POD is involved to
reduce the large dimension of the input snapshots and the GPR is applied as a
probabilistic response surface for the approximation of the parametric solution
manifold. The novelty of this work is in fact the creation of a computational
framework by combining model order reduction and probabilistic prediction, in
order to increase the accuracy of the state-of-the-art (data-based)
metodologies for a very limited amount of snapshots. 
This framework results then particularly suited for industrial setting, where
simulations can last days or even weeks, making unfeasible the computation of
large database of solutions.
Similar methods have be proposed in the last year to data-based modeling of
parametric problems, using interpolation to approximate the solution manifold.
Among all contributions in literature, we
cite~\cite{SalmoiraghiScardigliTelibRozza2018, demo2018shape,
tezzele2019marine} for a good overview. An alternative to interpolation or
regression for the online prediction is constituted by the Artificial Neural
Network (ANN), as described in~\cite{wang2019non}. Even if such approach
results very accurate, it requires larger datasets, resulting particularly
ineffective for limiting the number of full-order snapshots. Approaches that
exploit the Active Subspaces (AS) property for the reduction of the number of
required input snapshots are presented in~\cite{TezzeleDemoStabileMolaRozza2020MOR,DemoTezzeleRozza2019CRM}.

The paper is structured as follows: the POD and GPR methods are described more
in details in Section~\ref{sec:pod} and Section~\ref{sec:gpr}, respectively,
while in Section~\ref{sec:podgpr} we describe their integration into a single
framework. The numerical results are presented in Section~\ref{sec:results} and
they are subdivided as following: first we applied the POD-GPR framework to a
simply toy problem, a Stokes flow passing around a cylinder, then we move to an
industrial problem where turbulent biphase Navier Stokes flow is simulated
around a parametric cruise ship geometry. Finally, we summarize some
conclusion and further future perspectives in Section~\ref{sec:conclusion}.

\section{Proper orthogonal decomposition for dimensionality reduction}
\label{sec:pod}
POD is a technique used in the context of
PDEs aiming at reducing the
computational cost required to obtain their solution numerically. In this
section we introduce the reader to the general idea behind POD and give some
inshights on the mathematical theory supporting it.

We now begin with an introduction to the class of Reduced Basis (RB),
presenting POD in the following as a special case. Finally, we discuss a
modification of the standard POD formulation, which is equation free and it will
be the one used in this work.

\paragraph*{Reduced Basis.}
We begin with a formal definition of a system of discretized parametric PDEs in
the form: 

\begin{equation}
\text{find} \: u^\fulldim(\mupar) \in \fullspace \: \text{s.t.:} \:\:\:  a(u^\fulldim(\mupar),w;\mupar)=F(w;\mupar) \:\:\: \forall w \in \fullspace,
\label{truth_problem}
\end{equation}
where $\mupar \in \parspace \subset \mathbb{R}^p$ is the parameter, $u^\fulldim$ is the truth
solution to our problem, $a(\cdot, \cdot; \mupar)$ is the parametric bilinear
form and $\fullspace$ is a finite dimensional space of
dimension $\fulldim$.

We introduce the notion of the solution manifold, that is the set of all
possible solutions of our parametric problem under the variation of the
parameter:
\begin{equation}
\mathbb{M}^\fulldim=\{ u^\fulldim(\mupar) \:\: \mupar \in \parspace\}.
\end{equation}

The final goal of RB is to approximate any element of this manifold using
a low number of basis functions, or modes, $\{ \chi_i(x)\}_{i=1}^\reddim$ by setting
what we call the reduced basis. These $\reddim$ functions are globally defined over
the computational domain and are obtained using some pre-computed truth
solutions for particular parameter values. 

The \textit{reduced solution} $u_\reddim^\fulldim \approx u^\fulldim$ is
hence composed by a suitable linear combination of these modes:
\begin{equation}
u_\reddim^\fulldim(\mupar) = \sum_{i=1}^\reddim \xi_i(\mupar) \chi_i(x),
\end{equation}
and the reduced formulation of the problem (\ref{truth_problem}) becomes:
\begin{equation}
\text{find} \: u_\reddim^\fulldim(\mupar) \in \redspace \: \text{s.t.:} \:\:\: a(u_\reddim^\fulldim(\mupar),w;\mupar)=F(w;\mupar) \:\:\: \forall w \in \redspace,
\label{reduced_problem}
\end{equation}
where $\redspace = span(\{ \chi_i(x)\}_{i=1}^\reddim)$.

Hence, while the degrees of freedom associated with the truth solution
$\fulldim$ are typically high, the degrees of freedom associated with the RB
approximation are only $\reddim$, where $\reddim \ll \fulldim$.

\paragraph*{Proper Orthogonal Decomposition.}
POD refers to a particular RB which considers a hierarchical orthonormal
basis generated using energetic considerations, exploting them as global basis
function for a Galerkin framework~\cite{ballarin2015supremizer}. The
mathematical tool behind POD is Singular Value Decomposition (SVD), as we will
briefly discuss now.

Let us begin with some notation. We will call a possible outcome of our system
$\mathbf{y}$ a \textit{snapshot} and denote with $\fulldim$ its dimensionality
($\mathbf{y} \in \R^\fulldim$). We will then denote with $n$ the number of
snapshots collected, and gather them in the \textit{snapshots' matrix}
$\mathbf{Y} \in \mathbb{R}^{\fulldim\times n}$:
\begin{equation}
\mathbf{Y} = 
\left[
  \begin{array}{cccc}
    \vrule & \vrule & & \vrule\\
    \mathbf{y}_{1} &\mathbf{y}_{2} & \ldots & \mathbf{y}_{n} \\
    \vrule & \vrule & & \vrule 
  \end{array}
\right].
\end{equation}

We then call $\reddim$ the number of POD modes, i.e. the dimension of the reduced basis, and assume $\reddim \leq n<\fulldim$. 

A formal definition of the POD basis can be then given:
\begin{theorem}
\textbf{(POD basis)} Given $\mathbf{Y} \in \mathbb{R}^{\fulldim \times n}$, $\{\chi_i\}_{i=1}^\reddim$, for $\reddim \in \{1,..,n\}$ is the POD basis of $\mathbf{Y}$ if and only if it is a solution of:
\begin{equation}
\max_{\tilde{\psi_1},\tilde{\psi_2},..\tilde{\psi_\reddim}} \sum_{i=1}^\reddim \sum_{j=1}^n | <\mathbf{y}_j,\tilde{\psi_i} >_{\R^\fulldim} |^2 \:\:\:\:\: s.t. <\tilde{\psi_i},\tilde{\psi_i} >_{\R^\fulldim} = \delta_{i,j}, \:\:\: for \:\: 1 \leq i,j \leq \reddim.
\label{pod_basis}
\end{equation}
\end{theorem}

This can be read as: the POD basis is the one that maximizes the similarity (as
measured by the square of the scalar product) between the snapshots matrix and
its elements, under the constraint of orthonormality. In this sense, when we
obtain the $\reddim$-rank POD basis, we have the set of dimension $\reddim$ capable of
optimally express the variance in the snapshots.

The link between POD and SVD is stated in following theorem, that can be proven
using Lagrangian penalization techniques (see for example~\cite{volkwein-rom}):
\begin{theorem}
Given $\mathbf{Y} \in \mathbb{R}^{\fulldim \times n}$ of rank $d=n<m$, its
$\reddim$-rank POD basis is given by the set of the first $\reddim$ left singular
vectors $\{ \psi_i \}_{i=1}^\reddim$.
\end{theorem}

As a last remark on POD, we highlight the fact that, as already shown, to
compute the reduced basis we need a set of what can be called, using machine
learning nomenclature, training data, that is a set of truth solutions for
particular values of the parameters computed using the high fidelity solver.
This is usually the most expensive phase in the ROM pipeline, and the
applicability of ROMs techniques is often related to the size of the sampling
needed to obtain certain performances.

Connected to this last aspect, we mention here another important feature of
ROMs, that is the offline-online decomposition. By this, we mean that the use
of RB techniques can be split, as we have already seen, in two phases:
\begin{itemize}
\item an \textit{offline-phase}, in which a set of high fidelity solutions are
collected and the reduced basis is computed by combining them;
\item an \textit{online-phase}, in which we project the problem onto the
reduced space and the solutions for new parameters are computed in an efficient
manner.
\end{itemize}

The first phase requires typically more computational time, having to rely on
the high fidelity solver, and is usually executed in high performance computing
clusters using parallel programming techniques. The second phase, however, is
quite inexpensive from the computational point of view and could be done in
low-end terminal or even tablet/smartphones, widening the field of
applicability of this method.

\paragraph*{Data-driven Proper Orthogonal Decomposition.}
Data-driven POD is an alternative to standard POD that is more used in the
industrial setting because of its ease of use and its flexibility.  Within this
method, the first phase, the generation of the reduced basis, is performed in
the same way as the standard approach. A significant difference, however, is
that while in intrusive POD-Galerkin we need to rely upon an open-source
software to compute the truth solutions (since in the online phase we will need
to have access to the source code in order to project the equations). In the
data-driven approach, this is not necessary, and we can use commercial software
or even experimental data to train our model.  The reason why in data-driven
POD methods we have this freedom in the offline phase is that the online phase
is completely different from the standard POD.

We recall here that the assumption of RB-ROMs is that the truth solution of
our problem $u^\fulldim$ can be approximated by the reduced solution
$u_\reddim^\fulldim$ composed by linear combination of spatial modes
$\chi_i(x)$ multiplied by coefficients $\xi_i(\mupar)$, that is:
\begin{equation}
u^\fulldim(\mupar) \approx u_\reddim^\fulldim(\mupar)=\sum_{i=1}^\reddim \xi_i(\mupar) \chi_i(x).
\label{reco}
\end{equation}

Then we build an interpolation, or a regression, using high-fidelity data
coming from offline phase, approximating the function that associates the value
of the parameter $\mupar$ to the modal coefficients of the related solution $\{
\xi_i(\mu) \}_{i=1}^\reddim$.

We can use the interpolator to infer the value of the coefficients associated
with new parameters. The values of the coefficients are then used to
reconstruct the approximated truth solution using (\ref{reco}).

This approach is entirely data-driven and is independent both on the equations
and on the physics of the problem. This has its own advantages and disadvantages.
The ease of use and the complete freedom in the generation of the snapshots,
that are crucial in an industrial setting in which commercial software are
widely spread, correspond to a lower accuracy associated with the reduced
model. 

\section{Gaussian process regression model}
\label{sec:gpr}

We introduce in this section the Gaussian Process Regression (GPR), the supervised
learning technique we adopt to approximate the modal coefficients. We provide
in the following lines a summary of the method, in order to let the reader 
understand the entire pipeline, but avoiding to touch more details; for a
complete description of such framework, we
suggest~\cite{quinonero2005unifying}.

A Gaussian Process (GP) is a stochastic process whose finite-dimensional distributions are Gaussian. In supervised learning,
this process is usually adopted for regression manner, by constructing a
stochastic model from a set of input data capable to predict quantities of
interest. We define the set $\mathcal{D} = \{{(x_i, y_i)}\}$ for $i = 1, 2,
\dotsc, M$ as the set containing all the input-output pairs, where $x_i \in
\mathcal{P} \subset \R^m$ is the input parameter and $y_i \in \R$ is the
corresponding output.
We assume the output follows an unknown regression function $f : \mathcal{P}
\to \R$ we want to approximate (we avoid in this simple overview to deal with
noise), such that $y_i = f(x_i)$ for $i~=~1,2,\dotsc,M$. Thus, we define a GP
such that:
\begin{equation}
f(x) \sim \text{GP}(\mu(x), \mathbf{K}(x)),
\end{equation}
where $\mu(x)$ refers to the mean of the distribution and $\mathbf{K}(x)$ to
its covariance, defined as $\mathbf{K}_{ij} = \mathcal{K}(x_i, x_j)$. The
covariance function, also called {\it kernel}, $\mathcal{K}: \mathcal{P}
\times \mathcal{P} \to \R$ should be properly selected depending on the
problem, for this work we use the squared esponential:
\begin{equation}
\mathcal{K}(x_i, x_j) = \sigma^2 \exp\left(-\frac{\|x_i - x_j\|^2}{2l}\right).
\label{eq:se}
\end{equation}
The prior joint Gaussian distribution for the outputs $y$ results then
$$
\mathbf{y}|\mathbf{X}\sim~\mathcal{N}(0,\mathbf{K}_{\Xtrain\Xtrain}),
$$
where $\ytrain \in \R^M$ is the output vector and $\mathbf{X} \in
\R^{m\times M}$ is the matrix whose columns are the input parameters. We
specify that, for sake of compactness, we adopt a compress notation for the
covariance matrices, such that $\mathbf{K}_{\Xtrain\Xtrain} = (\mathcal{K}(x_i,
x_j))$ with $i, j = 1, \dotsc, M$, where $x_i$ and $x_j$ are the $i$-th and
$j$-th columns of $\Xtrain$.
We are now interested in predicting the output for new test input $\Xtest$
by exploiting the joint distribution based on the above prior, that is:
\begin{equation}
\left[
	\begin{array}{c}
		f(\Xtrain) \\
		f(\Xtest)
	\end{array}
\right] \sim \mathcal{N}
\left( \left[\begin{array}{c}\mu(\Xtrain) \\ \mu(\Xtest)\end{array}
\right],\,
\left[
	\begin{array}{ll}
	\mathbf{K}_{\Xtrain \Xtrain} & \mathbf{K}_{\Xtrain \Xtest}\\
	\mathbf{K}_{\Xtest \Xtrain} & \mathbf{K}_{\Xtest \Xtest} \\
	\end{array}
\right]
\right).
\end{equation}
Now, the ouputs $\ytest = f(\Xtest)$ can be expressed in probabilistic terms by simply using the conditional Gaussian distribution:
\begin{equation}
	\ytest | \Xtest, \Xtrain, \ytrain \sim \mathcal{N} (\mathbf{m},
	\mathbf{C}),\quad \mathbf{m} =
	\mathbf{K}_{\Xtest\Xtrain}\mathbf{K}_{\Xtrain\Xtrain}^{-1}\ytrain,
	\quad \mathbf{C} = \mathbf{K}_{\Xtest\Xtest} - \mathbf{K}_{\Xtest\Xtrain}\mathbf{K}_{\Xtrain\Xtrain}^{-1}\mathbf{K}_{\Xtrain\Xtest}.
\label{eq:gpr_pred}
\end{equation}
It is important to remark that, for a good result regarding the prediction,
we need to optimize the hyperparameters within the covariance function. In
this work, since a squared exponential kernel (Eq.~\ref{eq:se}) has been
adopted, we optimize the hyperparameters $\sigma$ and $l$, respectively the variance and the lenghtscale, by maximizing the likelihood
through a standard multi-start grandient-based optimization algorithm.

In this work, concerning the GPR, we use the Python open source
package called GPy~\cite{gpy2014}.

\section{Numerical results}
\label{sec:results}

\newcommand{\dd}{\ensuremath{\Omega}}
\newcommand{\p}{\ensuremath{p}}
\newcommand{\vv}{\ensuremath{u}}
\newcommand{\vecp}{\ensuremath{p_h}}
\newcommand{\vecv}{\ensuremath{u_h}}
\newcommand{\Vh}{\ensuremath{\mathbb{V}_h}}
\newcommand{\Ph}{\ensuremath{\mathbb{Q}_h}}

We present in this section the numerical results obtained by applying the
described method to two different computational fluid dynamics simulations. To
emphasize the versatility of the proposed framework with respect to the
high-fidelity solver, we firstly apply it to a finite element (FE) problem then
to a finite volume (FV) one. The coupling between the original model
and the reduced one is in fact based only on data, resulting in a complete
modular pipeline. In the next sections, we briefly introduce both the
full-order models before providing a discussion about the reduction accuracy, even
if in this case the original model has not a strong impact in the generation of
the reduced order space, in order to give to the reader the possibility to
reproduce the results.

\subsection{Parametric Stokes flow around cylinder}

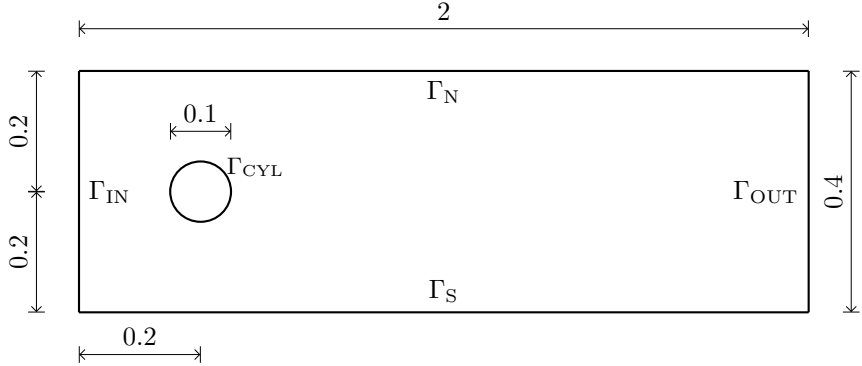
\begin{figure}[ht]
\centering
\begin{tikzpicture}[scale=.8,
measure/.style = {very thin,
        {Bar[width=2.2mm]Straight Barb[]}-%
        {Straight Barb[]Bar[width=2.2mm]},
            },
every label/.append style = {font=\footnotesize}
                    ]
\coordinate (a) at (0,0);
\coordinate (b) at (12,0);
\coordinate (c) at (12,4);
\coordinate (d) at (0,4);
\coordinate (e) at (2,2);
\coordinate (f) at (2,0);
\coordinate (g) at (0,2);
\coordinate[label={\small $\Gamma_\text{CYL}$}] (h) at (2.9, 2.1);
\draw[black, thick] (a) -- (b) node[midway, above] {$\Gamma_\text{S}$};
\draw[black, thick] (b) -- (c) node[midway, left] {$\Gamma_\text{OUT}$};
\draw[black, thick] (c) -- (d) node[midway, below] {$\Gamma_\text{N}$};
\draw[black, thick] (d) -- (a) node[midway, right] {$\Gamma_\text{IN}$};
\draw[black, thick] (e) circle (.5 cm);

\draw[measure]  ($(c)!7mm!-90:(d)$) to [sloped,"$2$"] ($(d)!7mm! 90:(c)$);
\draw[measure]  ($(a)!-7mm! 90:(f)$) to [sloped,"$0.2$"] ($(f)!-7mm!-90:(a)$);
\draw[measure]  ($(a)!7mm!90:(g)$) to [sloped,"$0.2$"] ($(g)!-7mm!90:(a)$);
\draw[measure]  ($(g)!7mm!90:(d)$) to [sloped,"$0.2$"] ($(d)!-7mm!90:(g)$);
\draw[measure]  ($(b)!-7mm!90:(c)$) to [sloped,"$0.4$"] ($(c)!7mm!90:(b)$);
\draw[measure]  (1.5, 3) to [sloped,"$0.1$"] (2.5, 3);
\end{tikzpicture}
\caption{The domain $\dd$ for the parametric Stokes flow simulation.}\label{fig:stokes_domain}
\end{figure}
The first numerical experiment where we use the POD-GPR framework is the
simulation of a parametric Stokes flow passing around a circular cylinder. We
define $\vv(\mupar):\dd\to\R^2$ and $\p(\mupar):\dd\to\R$ the velocity and the
pressure parametric fields, respectively, such that: 
\begin{equation}
\begin{array}{rll}
\mu(\mu_0)\Delta {\bf \vv}(\mupar) + \nabla \p(\mupar) &= 0\quad\quad\quad&\text{in}\,\Omega,\\
\text{div}\,\vv(\mupar) &= 0\quad\quad\quad&\text{in}\,\Omega,\\
\vv(\mupar) &= \mu_1\quad\quad\quad&\text{in}\,\Gamma_\text{IN},\\
\vv(\mupar) &= 0\quad\quad\quad&\text{in}\,\Gamma_\text{S} \cup \Gamma_\text{N}\cup \Gamma_\text{CYL}, \,\\
\p(\mupar) &= 0\quad\quad\quad&\text{in}\,\Gamma_\text{OUT}.\\
\end{array}
\label{eq:stokes}
\end{equation}
As we can see from Eq.~\ref{eq:stokes}, we consider two physical parameters
$\mupar = (\mu_0, \mu_1)\in \parspace \subset\R^2$ controlling the viscosity $\mu$ of the
fluid and the velocity of the fluid at the inlet boundary $\Gamma_\text{IN}$.
The 2-dimensional domain $\dd$ (sketched in Figure~\ref{fig:stokes_domain}) is
a rectangle with a circular hole, on which we impose on the physical walls
($\Gamma_\text{S} \cup \Gamma_\text{N}\cup \Gamma_\text{CYL}$) no-slip boundary
condition, on the inlet the parametric flow velocity $\mu_1$ and on the outlet
a null pressure.

Starting from Eq.~\ref{eq:stokes}, we can obtain the weak formulation of the
Stokes equations by multiplying these equations to the test functions $v \in
[H^1_0(\Omega)]^2 = V$ and $q \in L^2(\Omega) = Q$:
\begin{equation}
\begin{array}{rll}
a(\vv(\mupar), v; \mupar) + b(v, \p(\mupar); \mupar) &= 0\quad\quad\quad\forall v \in V,\\
b(\vv(\mupar), q; \mupar) &= 0\quad\quad\quad\forall q \in Q.\\
\end{array}
\end{equation}
where $a(\cdot, \cdot; \mupar)$ and $b(\cdot, \cdot; \mupar)$ are the parametric bilinear forms.

To solve the Stokes problem described above, we adopt a finite element (FE)
approach. We define the finite dimensional spaces $\Vh = (P_k(\mathcal{T}))^2
\cap V$ and $\Ph = P_l(\mathcal{T}) \cap Q$, where $P_i(\mathcal{T}) = \{v \in
C(\dd)\,:v|_\tau \in P_i, \forall \tau \in \mathcal{T}\}$ for $k > 1$ indicates
the piecewise polynomial space defined on the triangulation $\mathcal{T}$ of
the domain $\dd$. For this work, we have chosen a second order polynomial for
the velocity space and a first order one for the pressure space, i.e. $(P_2,
P_1)$. The use of the so-called Hood-Taylor scheme~\cite{taylor1973numerical}
is sufficient for ensuring the inf-sup condition; for the mathematical proof
of the stability of such scheme, we refer to~\cite{boffi2008finite}. Moreover,
the Reynolds number is very low ($\text{Re} < 1$) for all the combinations of
parameters, making possible to avoid additional stabilization to
numerically solve the problem.

The discretized parametric Stokes problem reads as follow:
\begin{equation}
\begin{array}{rll}
a(\vecv(\mupar), v_h; \mupar) + b(v_h, \p_h(\mupar); \mupar) &= 0\quad\quad\quad\forall v_h \in \Vh,\\
b(\vecv(\mupar), q_h; \mupar) &= 0\quad\quad\quad\forall q_h \in \Ph,\\
\end{array}
\end{equation}
where $\vecv(\mupar) \in \Vh$ and $\vecp(\mupar) \in \Ph$ are the finite
dimensional unknowns we want to compute.
The discretized problem has been solved using the
FEniCS framework~\cite{AlnaesBlechta2015a}.

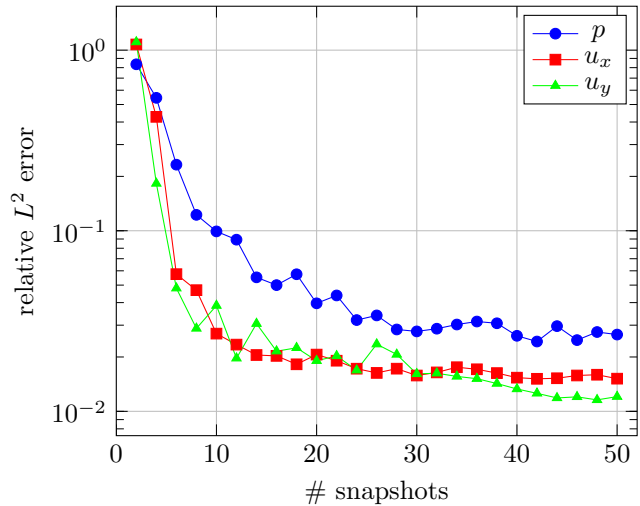
\begin{figure}
\centering
\begin{tikzpicture}

\begin{axis}[
tick pos=both,
xmin=0, xmax=52,
ymode=log,
xtick style={color=black},
ytick style={color=black},
xlabel={\# snapshots},
ylabel={relative $L^2$ error},
grid
]
\addplot [blue, mark=*]
table {%
2 0.8350522037428473
4 0.5444542154568947
6 0.2321842145610416
8 0.12245333771508393
10 0.0991743246202657
12 0.08924428114573797
14 0.0552136502644083
16 0.049995381057887804
18 0.05742737479533459
20 0.03964305360998087
22 0.043829701576101845
24 0.0320233190866012
26 0.033977447404885154
28 0.02840080279162784
30 0.027710208541722077
32 0.028681965545205045
34 0.030253742410927553
36 0.03142999403353226
38 0.03073916765149662
40 0.026198490508949995
42 0.024357293939024426
44 0.029636872661247628
46 0.024801862205239424
48 0.027493661347707133
50 0.026614299814170135
};
\addlegendentry{$p$}

\addplot [red, mark=square*]
table {%
2 1.0751235260089982
4 0.42731088244600085
6 0.05749245390179458
8 0.04687777728075716
10 0.026939589733566315
12 0.02341647136871812
14 0.020504693528810215
16 0.020279176894338814
18 0.018219203177425102
20 0.020619336407070234
22 0.019072808170106918
24 0.017226290646823874
26 0.01632232271987715
28 0.01723108069860893
30 0.015791051994561742
32 0.016429878509511726
34 0.01755202394264883
36 0.017098210047488786
38 0.01628751916577057
40 0.015362176989579309
42 0.015143279833557701
44 0.015249132494745563
46 0.015796680291792785
48 0.01594842452361992
50 0.015169699023775718
};
\addlegendentry{$u_x$}

\addplot [green, mark=triangle*]
table {%
2 1.1091774059727328
4 0.18274011209185434
6 0.04806719535318496
8 0.02878557410116909
10 0.03853608548643857
12 0.019678383511160226
14 0.030602018365674583
16 0.02153324309749566
18 0.022478856380663007
20 0.019061427121804763
22 0.020270467694966213
24 0.01688886163465804
26 0.023551224790325693
28 0.020653520723217512
30 0.016123748201219208
32 0.01625448319684291
34 0.015611334159424764
36 0.015164221389588556
38 0.014256348426184276
40 0.013299103350255668
42 0.012582665233954238
44 0.011864960415934853
46 0.012028723075658568
48 0.011574392516793413
50 0.012064001727256354
};
\addlegendentry{$u_y$}
\end{axis}

\end{tikzpicture}
\caption{Relative average error between the truth solutions and the solutions
obtained with the POD-GPR framework, evaluated in 20 test parameters, varying
the number of snapshots.}\label{fig:err_stokes}
\end{figure}
We initially create the high-fidelity database by computing the pressure and
velocity fields for 50 samples, randomly generated in the parameter space
$\parspace = [.1, 10]\times[1\mathrm{e}{-4}, 1\mathrm{e}{-6}]$. We then apply
the described method to these snapshots, using initially the POD algorithm to
reduce the dimensionality of the solutions, aka extract the modal
coefficients. 

\begin{figure}[htb!]
\includegraphics[trim=150 0 80 0,clip,width=\textwidth]{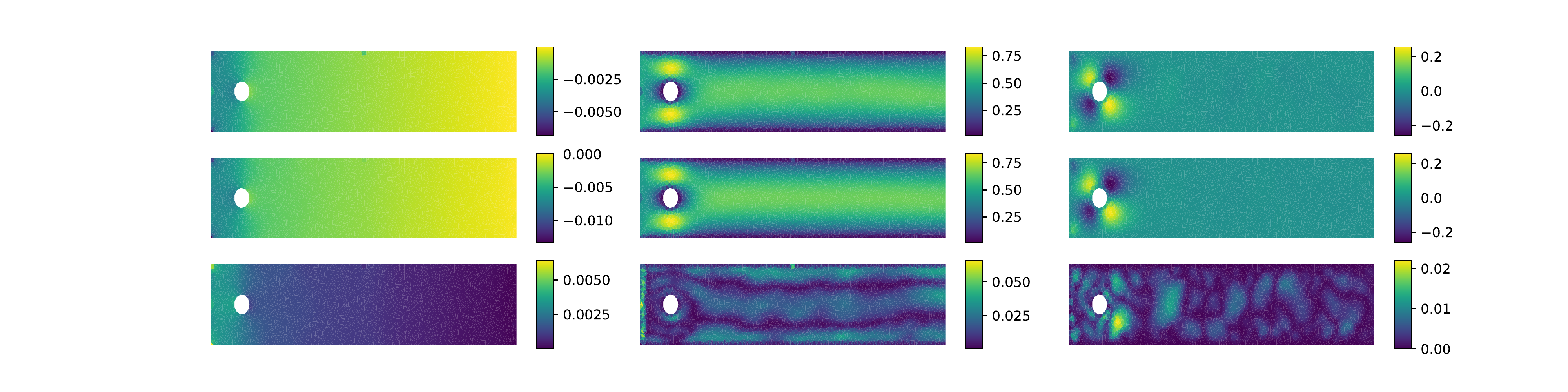}
\caption{Graphical visualization of the reconstructed fields for a test
parameter using only 3 input samples. The plots are organized as following:
the three columns represent, from left to right, pressure, velocity among $x$
direction and velocity among $y$ direction. The first row shows the truth
solutions, the second row the solutions approximated using the POD-GPR method
and the third one represents the absolute error between them.
}\label{fig:err_stokes2}
\end{figure}

In this case, we use 12 modes for the construction of the
reduced space and apply a GP --- we repeat that the used kernel is the
squared exponential one --- to the distribution of the modal coefficients.
To test the accuracy, we create a test dataset containing 20 random samples
selected in the parameter space, that of course are different to the input
database. For all these test parameters, we approximate the reduced solutions
and compute the realtive error between the reduced test solutions and the
truth ones. Figure~\ref{fig:err_stokes} reports the average relative error
for all the analized fields in relation to the number of input snapshots,
starting from only 2 snapshots.  We specify that, using less than 12
snapshots, the number of POD modes is limited to the number of snapshots.
As we can see, the error results very small with even few snapshots,
especially for both the velocity components that are beyond the 5\% with only
8 snapshots. The trend for the pressure error is less steep (we need 16
snapshots to reach an error of 5\%), but it is able to reach at the end a
precision of about 2.5\%. A graphical visualization of the test solutions
obtained using only 3 input samples is provided in
Figure~\ref{fig:err_stokes2}.
We underline the error becomes pretty much constant using many snapshots, due
to the Bayesian online phase, making the usage of GP during the online stage
profitable especially in the cases where the offline stage has a limited computational budget.

\subsection{Multiphase turbulent Navier-Stokes flow in an industrial parametric domain}

The second numerical experiment we set up to test the POD-GPR framework is
the simulation of the flow around the hull of a cruise ship. The reduction is
then employed inside a shape optimization pipeline, enabling the use of a
global optimization procedure on the design space. 

We now proceed to discuss the main ingredients of the optimization pipeline
in order to understand the role of the POD-GPR reduction. We avoid going into
many technical details and refer the reader to~\cite{mythesis, DemoOrtaliGustinRozzaLAvini2020BUMI} for a more
extensive presentation of all the optimization pipeline.

The first step of the optimization pipeline consists of the generation of the
parametric design space, where each possible shape is associated to a
particular value of a pre-defined parameter.  This is done by considering an
initial undeformed configuration and applying to it a technique called Free
Form Deformation (FFD) (for more details see~\cite{ffd2,ffd1,SalmoiraghiBallarinHeltaiRozza2016}).
FFD is a deformation technique in which the object to be deformed is
initially put inside a lattice of points; then some of these points are
moved, in a parametric way, and the space enclosed is deformed smoothly
according to these motions. 

\begin{figure}[ht!]
\centering
\includegraphics[width=.8\textwidth]{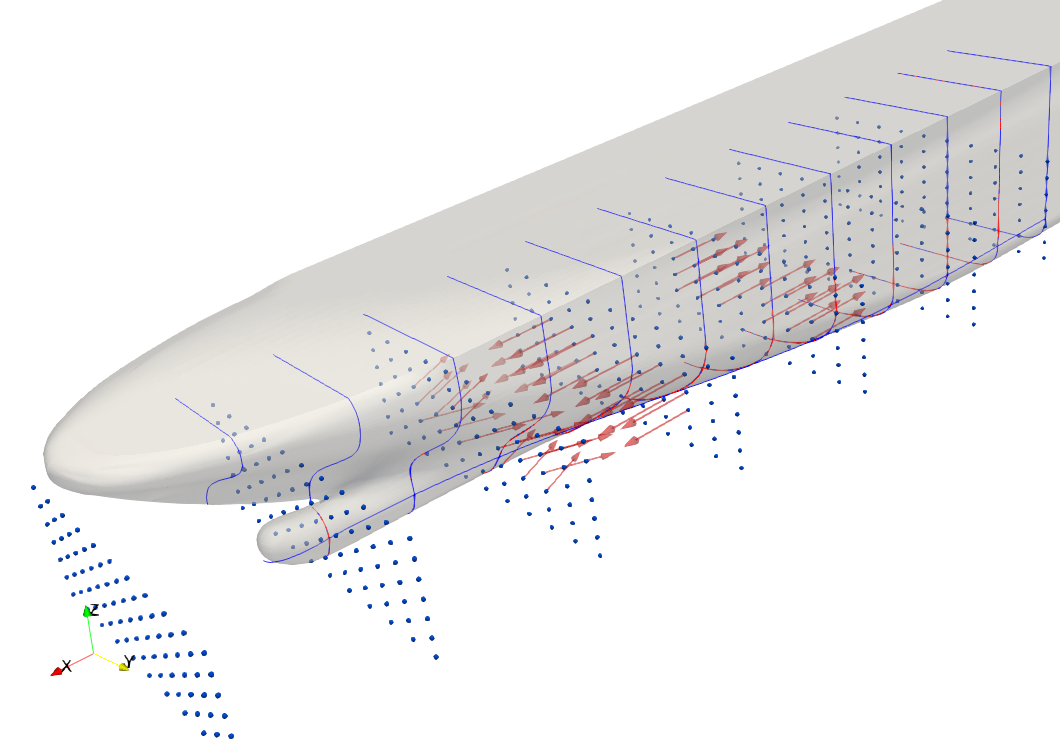}
\caption{View of the the lattice of points over the undeformed hull: the blue
dots are the FFD control points, while the red arrows represent the
displacements.}
\label{grid}
\end{figure}

In our case, we define the lattice shown in Figure~\ref{grid}. This lattice
is enclosing the bottom-frontal part of the boat, being the part we want to
deform, and is composed of a total of around 500 points, of which only a
small part will be displaced.

In particular we define a set of movements connected to 6 independent
parameters. The way in which these parameters are defined and the ranges they
take values from are related to the kind of shapes we want to consider and
the various constraints those shapes are subjected to. For example, the most
simple and straightforward constraint we imposed is that of volume, i.e we
impose that the hulls must not have a total volume that reaches too small
values. Another constraint we impose is that the deformations applied on the
hull needs to connect continuously and smoothly to the undeformed parts of
the boat. This latter constraint can be eimposed by keeping fixed some points
in the lattice.

The FFD step is performed using PyGEM, an open-source Python library
(see~\cite{pygem}). 

The next step consists of the evaluation of the performance of the newly
generated hull. This is done via numerical simulation of the flow of water
and air around it. In particular, we consider a 3D biphase (water and air)
incompressible flow. The physical model is based on the well known
Incompressible Navier-Stokes equations, with modellation of turbulence based on a RANS approach (for more details see~\cite{turbobook}), Finite Volume discretization (see~\cite{mouk}) and  using the Volume of Fluid (VOF) methodology (see~\cite{vof}) to take into account the biphase
nature of the fluid. 

More specifically, the equations considered are the following:
\begin{equation}
\begin{cases}
\frac{\partial(\rho u)}{\partial t} + \nabla \cdot (\rho u\otimes u ) + \nabla p  - \rho g -\nabla \cdot \nu \nabla u - \nabla \cdot R=0, \\
\nabla \cdot u =0, \\ 
\frac{\partial \alpha}{\partial t} + \nabla \cdot (u \alpha)=0, \\
\end{cases}
\label{vof}
\end{equation}

where $u$, $p$ denote the fluid velocity and pressure fields, $\rho$, $\nu$ and $g$ denote respectively density, dynamic viscosity and gravity acceleration, $R$ is the Reynolds stress tensor and $\alpha$ is the fraction of water.

The first equation in Eq.~\eqref{vof} expresses the momentum balance in the
infinitesimal volume and is essentially a reformulation of Newton'st $2^{nd}$
law of dynamics expressing the Eulerian fluid's acceleration ($\frac{\partial
u}{\partial t}$) as a result of convection ($(u \cdot \nabla)u$), a gradient of
pressure ($\frac{1}{\rho} \nabla p$) and diffusion ($\nabla \cdot \nu \nabla
u$). The second equation instead expresses the mass conservation, emerging as
solenoidality of the velocity field because of the incompressibility
hypothesis.

We note that both the velocity and pressure fields in this equations represent only the mean part, decoupled from the fluctuating part using a RANS approach. All the effects of the fluctuating part are contained in the Reynolds stress tensor $R$, which will be in then modeled as a function of the mean velocity field $u$ in order to close the equations.

The third equation is coming from the VOF modelling of the biphase flow. 

Volume of fluid is a free-surface modeling technique that allows to describe a multiphase fluid composed of two incompressible, isothermal immiscible fluids (water and air, in our case). For a more in-depth discussion on VOF, we refer to~\cite{vof}.

This method is based on a phase-fraction technique: a scalar variable, denoted by $\alpha$, is added to N-S equations, representing the fraction of water contained in the infinitesimal volume (or in the finite volume, when we will discretize the equations). This variable belongs to the interval $[0,1]$, where $\alpha=1$ represents a point in which water is present, $\alpha=0$ a point in which air is present and $\alpha \in (0,1)$ represents interface points. By its nature, $\alpha$ is a discontinuous variable, and so the discretization method must ensure that the interface is captured a small number of cells.

In addition to the momentum and mass balance equations, we have two equations where the density $\rho$ and the kinematic viscosity $\nu$ are defined using an algebraic formula expressing them as a convex combination of the corresponding water and air properties:
\begin{equation}
\rho = \alpha \rho_W + (1-\alpha) \rho_A,
\label{density}
\end{equation}
\begin{equation}
\nu = \alpha \nu_W + (1-\alpha) \nu_A.
\label{density}
\end{equation}

Finally, the third equation in the system \eqref{vof} is the equation required to close the system after the addition of the variable $\alpha$. In fact, this equation is simply a transport equation for the fraction of fluid in which only the convective term is considered.

As for the discretization and solution of the above equations, we used the implementation contained in the open-source CFD software OpenFOAM (see~\cite{openfoam}). 

The POD-GPR framework described above is then employed in order to reduce the computational time required for the resolution of the numerical problem associated to a new shape. 
Going into the details of the reduction, we consider as a snapshot, and hence as a solution, not the whole volumetric flow field but the the field of total resistance (viscous+pressure) over the hull.

The first step consists of the collection of the snapshots matrix, i.e. a set
of solutions of the Full Order Model for selected parameter values. In
particular we sampled the parameters space described above considering the
vertices ($2^6=64$ points) and other 36 random points inside the domain, for a
total of 100 snapsohts. Other 20 random snapshots have been generated for
purpose of testing.

The extraction of the modes and approximation step have then been conducted
using EZyRB~\cite{DemoTezzeleRozza2018EZyRB}, an open-source Python library for
data-driven model order reduction.
We now show some graphs reporting the relative $L^2$ error between FOM and
ROM: 
\begin{equation}
\bigg( \frac{\int_A(u^\mathcal{N}-u^\mathcal{N}_N)^2 dA}{\int_A u^{\mathcal{N}2}_N dA}\bigg)^{\frac{1}{2}}.
\label{err}
\end{equation}

\begin{figure}
\begin{center}
\begin{tikzpicture}

\definecolor{color1}{rgb}{0.75,0,0.75}
\definecolor{color2}{rgb}{0.75,0.75,0}
\definecolor{color0}{rgb}{0,0.75,0.75}

\begin{groupplot}[
	group style={
		group size=2 by 1,
		horizontal sep=40pt,
	},
	width=.48\textwidth,
	grid,
    yticklabel style = {font=\scriptsize,xshift=0.5ex},
    xticklabel style = {font=\scriptsize,yshift=0.5ex},
]
\nextgroupplot[
    legend style={%
        column sep = 8pt,
        legend columns = -1,
        legend to name = grouplegend,
    },
    ylabel=relative $L^2$ error,
    xlabel=\# modes,
]
\addplot [blue, mark=*]
table {%
1 0.0647670728215381
3 0.0594114426723079
5 0.0553456394401349
7 0.0532667059959067
9 0.0527579392911333
11 0.0523239237504958
15 0.0513630872338408
20 0.0506309633501817
25 0.0508040216608267
30 0.0510512903515633
40 0.0526220437864548
60 0.0553470812622952
80 0.0575928650625103
};
\addlegendentry{GP}
\addplot [green!50.0!black, mark=square*]
table {%
1 0.0588433775747572
3 0.058424609774319
5 0.0577508923376143
7 0.0582170598286952
9 0.060558049317481
11 0.0609708436444476
15 0.0602791913752857
20 0.0604120970259428
25 0.060586459662219
30 0.0605634747339095
40 0.0611118263030476
60 0.0613459129629952
80 0.0614499113214857
};
\addlegendentry{linear}
\addplot [red, mark=triangle*]
table {%
1 0.0703925136041706
3 0.0654666824881262
5 0.0669587562137952
7 0.071035751636373
9 0.0763485993900294
11 0.0790590236081722
15 0.0817540323440976
20 0.084119749957877
25 0.0856928197649325
30 0.0882591140120198
40 0.092130711517646
60 0.096455359991419
80 0.0973185320922635
};
\addlegendentry{RBF0}
\addplot [orange, mark=diamond*]
table {%
1 0.0684607927552032
3 0.0630048072293492
5 0.0584480600140945
7 0.0563286489732706
9 0.055679303004727
11 0.0551906974742444
15 0.0534088269404548
20 0.052594452801573
25 0.0524500856222952
30 0.0523535839017198
40 0.0522098258318905
60 0.0523267193629476
80 0.0527518969983762
};
\addlegendentry{RBF0.01}

\nextgroupplot[
    xlabel=\# snapshots]
\addplot [blue, mark=*]
table {%
10 0.0641873803711008
15 0.0597589714254302
20 0.0582636993309825
30 0.0566366437175127
40 0.0553956167437841
60 0.0543279674642333
80 0.0535621058816762
};
\addplot [red, mark=triangle*]
table {%
10 0.0983259629012611
15 0.0865196690764611
20 0.0874302637727341
30 0.0837367872824286
40 0.0804280499747817
60 0.0779689702422913
80 0.0809664291099436
};
\addplot [orange, mark=diamond*]
table {%
10 0.0947043352924333
15 0.0770179270836516
20 0.0717440338652698
30 0.0637124012678087
40 0.0604545233550532
60 0.0567641143177484
80 0.0549787989854238
};
\end{groupplot}
\node at ($(group c2r1) + (-3.4cm,-4.1cm)$) {\ref{grouplegend}}; 
\end{tikzpicture}
\end{center}
\caption{Average relative $L^2$ error between the FOM and the ROM as a
function of the number of modes, keeping fixed the number of snapshots to 80,
(left) and the number of snapshots used, keeping fixed the number of modes to 20 (right). }
\label{err2}
\end{figure}
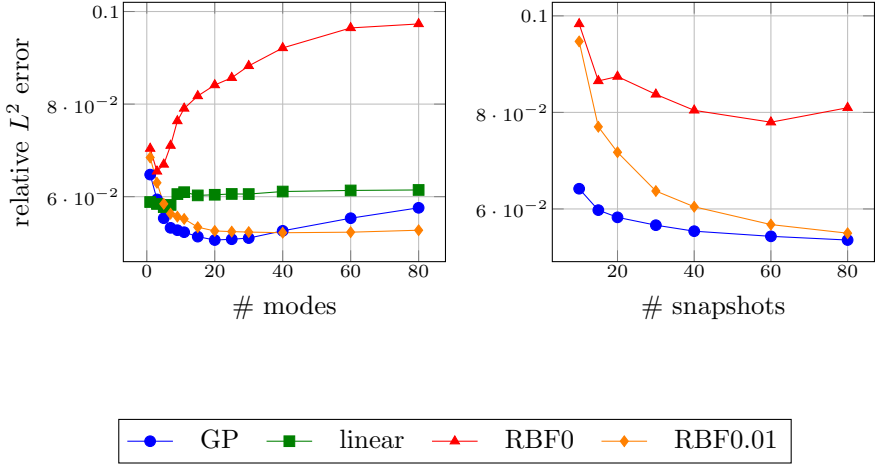

In particular, in Figure~\ref{err2} (left), the error is expressed as a
function of the number of modes used and as a function of the number of
snapshots used.  In both the figures we report in different colors different
regression techniques used in the approximation of the modal coefficients. We
consider the GPR approach, the classical $n$-dimensional linear
interpolation and two approximation techniques based on radial basis
functions (RBF), implemented using the Python package Scipy (see~\cite{2020SciPy-NMeth}). The former two methods differ one from the other for
the value of one particular parameter that defines the method, called
"smoothness", which quantifies how close the approximation technique is to
interpolation, where a value of 0 correspond to a plain interpolation while
higher values correspond to more smooth approximations. In particular, we
considered values of the smoothness equal to $0$ and $0.01$.

We can see from Figure~\ref{err2} that, using the GPR for the
regression of the modal coefficients, we are able to use a very low number of
snapshots obtaining an error which is comparable to the performances of
different regression techniques for higher number of snapshots. This is
really interesting since, as argued in the previous section, the offline
phase, i.e. the generation of the snapshots, is indeed the most expensive
from the computational point of view. However, when the number of snapshots
increases, the performances of GPR and the other methods tend to converge to a
same value.

The final step of the shape optimization pipeline consists on the exploration
of the design space in order to find the optimum. To perform this, we used a
evolutionary optimization algorithm, implemented inside DEAP, an open-source
Python library for evolutionary optimization (see~\cite{deap}). 

From the point of view of the speed-up in the computational time obtained with
the use of the POD-GPR framework, we mention that while in order to run the
above mentioned optimization procedure using the FOM we would have needed
around 140 weeks of CPU time, the use of the ROM approach restricted this time
to around 5 weeks. Of these 5 weeks, the biggest partition is associated to the
offline phase, i.e. the generation of the 100 training simulations, being the
time associated to the evaluation of the online phase almost negligible. It is
then clear that, in light of the considerations made regarding the error as a
funcion of the number of snapshots for the POD-GPR framework, we could have
augmented this speed-up of at least a factor 10 by generating fewer snapshots.

\begin{figure}
\centering
\includegraphics[width=.95\textwidth]{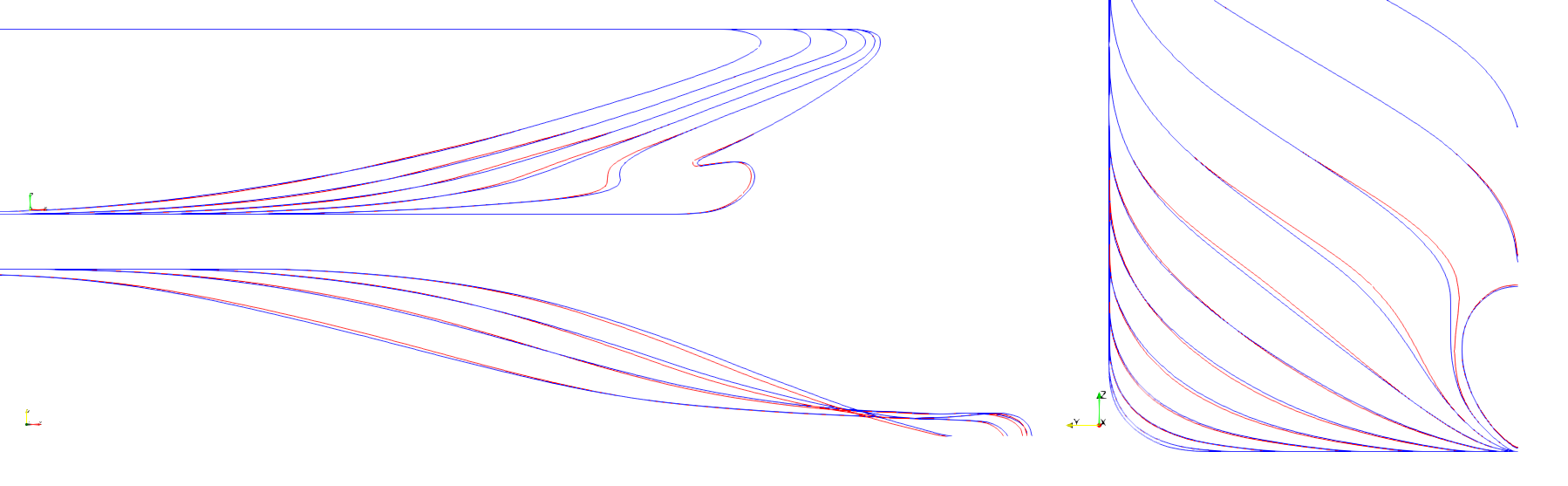}
\caption{View of the $x$ (right), $y$ (upper left) and $z$ (lower left) sections for the optimal hull shape; in blue the undeformed hull, in red the deformed hull.}
\label{opt_sez}
\end{figure}

Finally we report, in Figure~\ref{opt_sez}  the optimal shape found, which
correspond to a reduction of drag of $3\%$ from the initial hull. The figure
represents sections of the hull on the longitudinal and transversal plane, in
blue for the undeformed hull and red for the optimal deformed hull.

\section{Conclusions and perspectives}
\label{sec:conclusion}

In this contribution, we applied the POD-GPR framework to two completely
different parametric problems. To emphasize the equation-free nature of the
algorithm, we demonstrate the applicability to a Stokes problem discretized
using finite element method, then to a turbulent multiphase industrial problem
discretized using finite volume method. In both cases, it emerges that the GPR
method for the online prediction of the modal coefficients (for untested
parameters) is able to make a prediction with a relative small error by using very few input
snapshots. The advantage of such method is then maximazed when the offline
phase results quite expensive due to the high-fidelity model. 

The adoption of this method open several possibilities for future developments, by
using the GPR not only for the prediction of new modal coefficients, but also
for the quantification of the error, or for a greedy approach for new snapshots
selection by exploiting the variance of the Gaussian distribution. Additionaly, the use of different
kernels can be investigated.

\section*{Acknowledgements}
This work was partially funded by European Union Funding for
Research and Innovation --- Horizon 2020 Program --- in the framework
of European Research Council Executive Agency: H2020 ERC CoG 2015
AROMA-CFD project 681447 ``Advanced Reduced Order Methods with
Applications in Computational Fluid Dynamics'' P.I. Gianluigi Rozza.

\bibliographystyle{abbrv}
\bibliography{gpr-biblio}


\end{document}